\documentclass{article}
\usepackage{amsmath}
\usepackage{amssymb}
\usepackage{graphicx} 
\usepackage{caption,subfigure,subcaption}   
\usepackage{float}
\usepackage{mathtools}
\usepackage[usenames, table]{xcolor}
\usepackage{hyperref}

\usepackage{tabularx}      
\usepackage{array}         
\usepackage{booktabs}      
\usepackage{ragged2e}      
\definecolor{headerblue}{RGB}{220, 230, 242}
\definecolor{rowgray}{RGB}{245, 245, 245}
\definecolor{highlightyellow}{RGB}{255, 242, 204}
\newcolumntype{L}{>{\RaggedRight\arraybackslash}X}
\newcolumntype{C}{>{\centering\arraybackslash}X}

\definecolor{darkgreen}{rgb}{0,0.5,0}

\newcommand*{\dd}{\mathop{}\!\mathrm{d}}

\title{A weak order 2 Runge-Kutta method for It\^o stochastic delay differential equations}
\author{Alessia Andò, Dimitri Breda and Faraz William}

\date{}

\begin{document}

\maketitle
\begin{center}
CDLab -- Computational Dynamics Laboratory \\
Department of Mathematics, Computer Science, and Physics\\
University of Udine, Italy
\end{center}
\begin{abstract}

We present a Runge–Kutta method of weak order 2 for the numerical time integration of stochastic delay differential equations. This scheme extends the class of second order Runge–Kutta methods introduced by A. R{\"o}{\ss}ler in [SIAM J. Numer. Anal., 47(3):1713–1738, 2009] for stochastic ordinary differential equations. The proposed integrator is applicable to equations with discrete commensurable delays and is particularly efficient for problems involving multiple noise terms. Experimental confirmation of the weak order 2 is provided and MATLAB codes are freely available.

\end{abstract}



\section{Introduction}
\label{sec:intro}
Stochastic Delay Differential Equations (SDDEs) hold significant importance in the advanced modeling of complex systems, ranging from population dynamics and epidemics outbreaks in the natural sciences \cite{Mohammed1998Memory, Rihan2021SDDE} to supply chains for sustainable production systems \cite{supplychain2020}. They can be seen as a generalization of Stochastic Ordinary Differential Equations (SODEs) to portray phenomena in which, additionally to stochastic perturbations, the time elapsed between a cause and its effect(s) is taken into account. Specific instances are, e.g., the logistic equation, the predator-prey system and the brand goodwill supply chain described, respectively, in  \cite{Wu2012StochasticDelayLogistic}, \cite{predatorprey2015} and \cite{supplychain2020}.

Analytical expressions for the solutions of corresponding initial value problems (IVPs) of SODEs or SDDEs are rarely present. Therefore, numerical methods for their time integration are required to understand relevant behaviors and dynamics, but they pose in general nontrivial challenges arising from either stochasticity and memory effects, or their combination.

On the one hand, several numerical schemes are available for the time integration of SODEs, see, e.g., {\color{black}\cite{dambrosio2023numerical, higham2021introduction, KloedenPlaten1999, Platen_1999}} and the references therein. {\color{black}Particular attention has been devoted to their stability, in the mean-square and asymptotic sense \cite{ higham2000mean,saito1996stability} and, more recently, in the nonlinear framework through mean-square contractivity \cite{dambrosio2021contractivity}. The development of stabilized numerical integrators has also been investigated, for instance through the use of TASE operators to improve the stability properties of stochastic Runge--Kutta methods \cite{contemontano2026}. Furthermore, attention has been devoted to the preservation of structural properties of the underlying dynamics, such as the long-term conservation of invariants \cite{ChenCohenDAmbrosioLang2020,Debrabant2022,dambrosio2023longterm,Luesink2024}}. More recently, several higher order schemes have also been developed, including \cite{Debrabant2008,KloedenPlaten1999,Komori1997RootedTreeROW,Rossler2004StratonovichRK,Rossler2006RootedTreeSRK,Rossler2006scalarnoise,rossler2009,Tocino2002}. In particular, the class of Runge-Kutta (RK) methods presented by R{\"o}{\ss}ler \cite{rossler2009} for It\^o SODEs is specifically advantageous for its increased efficiency in dealing with multidimensional noise\footnote{\label{noisefootnote} For an $m$-dimensional noise only $2m-1$ random variables have to be generated rather than the $m(m+1)/2$ ones required by preceding methods like \cite{Debrabant2008,KloedenPlaten1999, Tocino2002}.}. Indeed, unlike all the other methods mentioned above, the RK schemes in this class share the distinguishing feature of a number of stages independent of the dimension of the noise. These methods are developed for It\^o SODEs, as the relevant order conditions are obtained by considering a sufficient number of terms in the underlying It\^o-Taylor expansion following It\^o's formula, which does not hold for Stratonovich SODEs.

On the other hand, the literature concerning the time integration of SDDEs is much sparser. The classical paper by Baker and Buckwar \cite{BakerBuckwar2000} is among the first works, where an extension of the celebrated Euler-Maruyama scheme from SODEs to SDDEs is presented along with a rigorous proof of convergence. Much more recently, the work by Zhou et. al  \cite{Mao2026} introduced a truncated version of the Euler-Maruyama framework for stochastic functional differential equations with infinite delay.
In between, a collection of papers appeared devoted to several aspects and improvements. Just as a couple of instances, \cite{MaoSabanis2003SDDE} presents a variant of Euler-Maruyama for SDDEs with variable time delay showing convergence under local Lipschitz conditions; \cite{WuMaoChen2009CIRDelay} shows the existence of strong solutions of a model with non-linear and non-Lipschitz diffusion coefficients, showing strong convergence of the relevant Euler-Maruyama approximations (but see also \cite{Buckwar2000SDDEIntro,kuchlerplaten,Mao1994VariableDelayPartII}). 

Most of the above literature concerning It\^o SDDEs deal with variants of the Euler-Maruyama scheme. However, the latter is limited to weak order 1. While, e.g., \cite{Li2024WeakTaylorSDDE} presents a method based on a It\^o Taylor symmetrical scheme that attains weak order 2, to the best of the authors' knowledge there are  no derivative-free methods for It\^o SDDEs that attain weak order 2. As part of the motivation of the present work we thus aim at extending from SODEs to SDDEs the class of RK methods developed in \cite{rossler2009}.

To complete a global picture of our motivation, we outline some directions that follow, even though a full treatment is target of ongoing and future research, and therefore beyond the scope of this paper. Indeed, an extension from deterministic DDEs to SDDEs of the class of Functional RK (FRK) methods originally introduced in \cite{MasetTorelliVermiglio} (for a recent self-contained and complete account see \cite{MasetDDE}) is work in progress of the authors. FRK schemes share nontrivial advantages with respect to more traditional time integrators for DDEs (see, e.g., \cite{BellenMasetZennaroGuglielmi2009,BellenZennaro2003}), which we believe they transfer to SDDEs as well. Firstly, they can tackle any functional dependence of the right-hand side on the past as they provide a numerical approximation of the solution continuously at any point of the current time step; this would allow one to treat more realistic models of the phenomena of interest, going beyond the reference case of a single, constant, discrete delay. Secondly, their formulation is such that explicit methods remain explicit in case of so-called overlapping, i.e., when the step size is larger than the minimum delay thus requiring the availability of solution values not yet computed. Let us note that advancing the time integration with explicit schemes and large step sizes is particularly effective when accuracy and stability requirements are satisfied, assuming that all breaking points of the desired order are anyhow included in the underlying mesh \cite{BellenMasetZennaroGuglielmi2009}. This would allow for an effective treatment of the more realistic models mentioned above, with due attention to SDDEs with multidimensional Wiener processes. In any case these advantages are more evident when convergence of high order is concerned: indeed, the functional explicit Euler method basically consists in using the classical explicit Euler method with any step size ranging from zero to the current step length. This explains our interest first in higher order RK schemes for SDDEs. In particular, the contents of the present manuscript regard a practical extension of the class \cite{rossler2009} to SDDEs with constant discrete delays, a relevant MATLAB implementation (which, to our knowledge, lacks for SODEs, too) freely available at \url{https://cdlab.uniud.it/software#w_ord_2_RK_SDDEs} and a thorough experimental validation on several instances of SDDEs confirming the preservation of the weak order 2.

\bigskip
The contents of the paper are organized as follows. In Section \ref{sec:rosslerSODE} we summarize from \cite{rossler2009} the class of RK methods of R{\"o}{\ss}ler for SODEs. In Section \ref{sec:rosslerSDDE} we illustrate the proposed extension to SDDEs. In Section \ref{sec:tests} we collect a series of numerical experiments providing validation of the weak order 2: they concern a linear SDDE, a nonlinear SDDE, a nonlinear system of SDDEs with multiple noises and, finally, a linear SDDE with multiple delays modeling a supply chain. Concluding remarks and future developments are discussed in Section \ref{sec:conclusions}.

\section{R{\"o}{\ss}ler's method for SODEs}
\label{sec:rosslerSODE}
For $d$ and $m$ positive integers and $T$ a positive real consider the SODE
\begin{equation}\label{SODE}
\dd y(t)=f(t,y(t))\dd t+\sum_{k=1}^{m}g_{k}(t,y(t))\dd W_{k}(t),\quad t\in[0,T],
\end{equation}
where $f\colon[0,T]\times\mathbb{R}^{d}\rightarrow\mathbb{R}^{d}$ represents the drift term and 
$g_{k}\colon [0,T]\times\mathbb{R}^{d}\rightarrow\mathbb{R}^{d}$ represents the diffusion coefficient associated to the $k$-th Brownian motion $W_{k}$, $k=1,\ldots,m$\footnote{Note that the stochastic part of \eqref{SODE} can be equivalently and more compactly described as $g(t,y(t))\dd W(t)$ for $g\colon [0,T]\times\mathbb{R}^{d}\rightarrow\mathbb{R}^{d\times m}$ and $W\in\mathbb{R}^{m}$, yet we keep the more explicit form in \eqref{SODE} in accordance with \cite{rossler2009}.}. Recall that a standard Brownian motion $t\mapsto W(t)\in\mathbb{R}$ is a Wiener process for which the random variable $W(t)$ satisfies
$W(0)=0$; for any $0\leq s<t\leq T$, $W(t)-W(s)$ is $\mathcal{N}(0,t-s)$ (or, equivalently, $\sqrt{t-s}\cdot\mathcal{N}(0,1)$); for any $0\leq s<t\leq u<v\leq T$, $W(t)-W(s)$ and $W(v)-W(u)$ are independent.

We assume throughout that $f$ and $g_{k}$, $k=1,\ldots,m$, are Lipschitz continuous and satisfy a standard linear growth condition with respect to their second argument. This guarantees existence and uniqueness of a strong solution of the Initial Value Problem (IVP) for \eqref{SODE} obtained by assigning a random variable $\phi\in\mathbb{R}^{d}$ as initial value, i.e., $y(0)=\phi$ (see, e.g., \cite[Chapter 4]{KloedenPlaten1999}). {\color{black}
Additionally, as required in \cite[Theorem 5.1]{rossler2009} to achieve weak convergence of order 2, we assume that
$f,g_k\in\mathcal{C}_{P}^{3,6}(\mathcal{I}\times\mathbb{R}^{d},\mathbb{R})$\footnote{{\color{black}$g\in\mathcal{C}_{P}^{l}(\mathbb{R}^{d},\mathbb{R})$ if its partial derivatives exhibit polynomial growth, i.e., there exist constants $C>0$ and $r\in\mathbb{N}$ such that
$
|\partial_x^i g(x)|\leq C(1+\|x\|^{2r})
$
for all $x\in\mathbb{R}^{d}$ and all $i\leq l$.
For $\mathcal{I}=[0,T]$, $g\in\mathcal{C}^{k,l}(\mathcal{I}\times\mathbb{R}^{d},\mathbb{R})
$ belongs to $\mathcal{C}_{P}^{k,l}(\mathcal{I}\times\mathbb{R}^{d},\mathbb{R})$ if
$g(t,\cdot)\in\mathcal{C}_{P}^{l}(\mathbb{R}^{d},\mathbb{R})$, with uniform polynomial growth with respect to $t\in\mathcal{I}$.}} for $k=1,\ldots,m$.
}

\bigskip
We summarize next the main ingredients of the class of RK methods developed in \cite{rossler2009}, proposed for It\^o SODEs with multiple noises such as \eqref{SODE} to overcome the drawbacks of previously introduced RK methods \cite{Debrabant2008,KloedenPlaten1999,Tocino2002} as explained in Section \ref{sec:intro}.

In what follows, let $y$ be the solution of \eqref{SODE} with $y(0)=\phi$ for a given $\phi$. For $N$ a positive integer, consider a constant step size $h\coloneqq T/N$ and introduce the uniform mesh $t_{n}\coloneqq nh$, $n=0,1,\ldots,N$, discretizing the time window $[0,T]$. Correspondingly, let $Y_{n}$ be the approximation of $y(t_{n})$ obtained for $n=0,1,\dots,N-1$ starting from $Y_{0}=y(0)$.
In order to construct such approximations, we first need to introduce the following independent random variables. First, for $k=1,\ldots,m$ let us define $\hat{I}_{k}$ with probability distribution
\begin{equation}\label{Ihat}
\left\{\setlength\arraycolsep{0.1em}\begin{array}{l}
\mathbb{P}\left(\hat{I}_{k}=\sqrt{3h}\right)=\frac{1}{6}\\[2mm]
\mathbb{P}\left(\hat{I}_{k}=-\sqrt{3h}\right)=\frac{1}{6}\\[2mm]
\mathbb{P}\left(\hat{I}_{k}=0\right)=\frac{2}{3},
\end{array}
\right.
\end{equation}
whose moments for $q\geq1$ integer have expected values
\begin{equation*}
\mathbb{E}\left(\hat{I}_{k}^{q}\right)=
\left\{\setlength\arraycolsep{0.1em}\begin{array}{ll}
0&\text{ for }q\in\{1,3,5\}\\[2mm]
(q-1)h^{q/2}&\text{ for }q\in\{2,4\}\\[2mm]
\mathcal{O}(h^{q/2})&\text{ for }q\geq6.
\end{array}
\right.
\end{equation*}
Second, for $k=1,\ldots,m$ let us define $\tilde{I}_{k}$ with probability distribution
\begin{equation}\label{Itilde}
\left\{\setlength\arraycolsep{0.1em}\begin{array}{l}
\mathbb{P}\left(\tilde{I}_{k}= \sqrt{h}\right)=\frac{1}{2}\\[2mm]
\mathbb{P}\left(\tilde{I}_{k}=- \sqrt{h}\right)=\frac{1}{2}
\end{array}
\right.
\end{equation}
whose moments for $q\geq1$ integer have expected values
\begin{equation*}
\mathbb{E}\left(\tilde{I}_{k}^{q}\right)=
\left\{\setlength\arraycolsep{0.1em}\begin{array}{ll}
0&\text{ for }q\in\{1,3\}\\[2mm]
h&\text{ for }q=2\\[2mm]
\mathcal{O}(h^{q/2})&\text{ for }q\geq4.
\end{array}
\right.
\end{equation*}
Finally, for $1\leq k,l\leq m$ let us define
\begin{equation}\label{Ihat2}
\hat{I}_{k,l}\coloneqq
\left\{\setlength\arraycolsep{0.1em}\begin{array}{ll}
\frac{1}{2}(\hat{I}_{k}\hat{I}_{l}-\sqrt{h}\tilde{I}_{k})&\text{ for }k<l\\[2mm]
\frac{1}{2}(\hat{I}_{k} {\color{black}\hat{I}_{l}}+\sqrt{h}\tilde{I}_{l})&\text{ for }k>l\\[2mm]
\frac{1}{2}(\hat{I}_{k}^2-h)&\text{ for }k=l.
\end{array}
\right.
\end{equation}
It is straightforward to check that the first five moments of $\hat{I}_{k}$ are the same as those of a true Brownian increment, i.e., $\mathbb{E}(\hat{I}_{k})=\mathbb{E}(\hat{I}^{3}_{k})=\mathbb{E}(\hat{I}^{5}_{k})=0$, $\mathbb{E}(\hat{I}^{2}_{k})=h$ and $\mathbb{E}(\hat{I}^{4}_{k})=3h^2.$ Moreover, $\tilde{I}_{k}$ matches the first three moments, i.e., {\color{black}$\mathbb{E}(\tilde{I}_{k})=\mathbb{E}(\tilde{I}^{3}_{k})=0$ }and {\color{black}$\mathbb{E}(\tilde{{I}}^{2}_{k})=h$}. Observe that these are the only moments playing a role in the weak error expansion up to order 2 (see, e.g., \cite[Section 14.2]{KloedenPlaten1999}).

Based on the above random variables, the numerical solution via R{\"o}{\ss}ler's method reads 
\begin{equation}\label{rosslerSODE:0}
\setlength\arraycolsep{0.1em}\begin{array}{rcl}
Y_{n+1}=Y_{n}&+&\displaystyle h\sum_{i=1}^{s}b_{i}f\left(t_{n}+c_{i}^{(0)}h,H_{n,i}^{(0)}\right)\\[4mm]
&+&\displaystyle\sum_{i=1}^{s}\sum_{k=1}^{m}\beta_{i}^{(1)}g_{k}\left(t_{n}+c_{i}^{(1)}h,H_{n,i}^{(k)}\right)\hat{I}_{k}\\[4mm]
&+&\displaystyle \sum_{i=1}^{s}\sum_{k=1}^{m}\beta_{i}^{(2)}g_{k}\left(t_{n}+c_{i}^{(1)}h,H_{n,i}^{(k)}\right)\frac{\hat{I}_{k,k}}{\sqrt{h}}\\[4mm]
&+&\displaystyle\sum_{i=1}^{s}\sum_{k=1}^{m}\beta_{i}^{(3)}g_{k}\left(t_{n}+c_{i}^{(2)}h,\hat{H}_{n,i}^{(k)}\right)\hat{I}_{k}\\[4mm]
&+&\displaystyle\sum_{i=1}^{s}\sum_{k=1}^{m}\beta_{i}^{(4)}g_{k}\left(t_{n}+c_{i}^{(2)}h,\hat{H}_{n,i}^{(k)}\right)\sqrt h, 
\end{array}
\end{equation}
where, for $i=1,\dots,s$ and $k=1,\dots ,m$,
\begin{equation}\label{rosslerSODE:1}
\setlength\arraycolsep{0.1em}\begin{array}{rcl}
H_{n,i}^{(0)}= Y_{n}&+&\displaystyle h\sum_{j=1}^{s}A_{i,j}^{(0)}f\left(t_{n}+c_{j}^{(0)}h,H_{n,j}^{(0)}\right)\\[4mm]
&+&\displaystyle\sum_{j=1}^{s}\sum_{l=1}^{m}B_{i,j}^{(0)}g_{l}\left(t_{n}+c_{j}^{(1)}h,H_{n,j}^{(l)}\right)\hat{I}_{l},
\end{array}
\end{equation}
\begin{equation}\label{rosslerSODE:2}
\setlength\arraycolsep{0.1em}\begin{array}{rcl}
H_{n,i}^{(k)}=Y_{n}&+&\displaystyle h\sum_{j=1}^{s}A_{i,j}^{(1)}f\left(t_{n}+c_{j}^{(0)}h,H_{n,j}^{(0)}\right)\\[4mm]
&+&\displaystyle\sum_{j=1}^{s}B_{i,j}^{(1)}g_{k}\left(t_{n}+c_{j}^{(1)}h,H_{n,j}^{(k)}\right)\sqrt h, \end{array}
\end{equation}
\begin{equation}\label{rosslerSODE:3}
\setlength\arraycolsep{0.1em}\begin{array}{rcl}
\hat{H}_{n,i}^{(k)}=Y_{n}&+&\displaystyle h\sum_{j=1}^{s}A_{i,j}^{(2)}f\left(t_{n}+c_{j}^{(0)}h,H_{n,j}^{(0)}\right)\\[4mm]
&+&\displaystyle\sum_{j=1}^{s}\sum_{\substack{l=1 \\ l\ne k}}^{m}B_{i,j}^{(2)}g_{l}\left(t_{n}+c_{j}^{(1)}h,H_{n,j}^{(l)}\right)\frac{\hat{I}_{k,l}}{\sqrt{h}}.
\end{array}
\end{equation}
Above, $s$ is the number of stages, $H_{n,i}^{(0)}$ and $H_{n,i}^{(k)},\hat{H}_{n,i}^{(k)}$ for $k=1,\ldots,m$ are the stage values, $b_{i}$ and $\beta_{i}^{(\ast)}$ are the weights, $c_{i}^{(\ast)}$ the abscissae and, finally, $A_{i,j}^{(\ast)}$ and $B_{i,j}^{(\ast)}$ the coefficients of the method. Correspondingly, we refer to the extended Butcher's tableau
\begin{equation*}
\renewcommand{\arraystretch}{1.3} 
\begin{array}{c|c|c|c}
c^{(0)} & A^{(0)} & B^{(0)}&\\
\hline
c^{(1)} & A^{(1)} & B^{(1)}&\\
\hline
c^{(2)} & A^{(2)} & B^{(2)}&\\
\hline
& b^{T} & \beta^{(1)T} & \beta^{(2)T}\\
\hline
&& \beta^{(3)T} & \beta^{(4)T}       
\end{array}
\end{equation*}
for (column) vectors $c^{(\ast)}$, $b$ and $\beta^{(\ast)}$ in $\mathbb{R}^{s}$ and matrices $A^{(\ast)}$ and $B^{(\ast)}$ in $\mathbb{R}^{s\times s}$.

In \cite{rossler2009} it is shown that the RK method \eqref{rosslerSODE:0}--\eqref{rosslerSODE:3} attains weak order 2 assuming to fulfill the order conditions described in Theorem 5.1 therein (to which we refer the reader for further details). We recall that the weak error is defined as
\begin{equation}\label{weakerror}
e_{\text{weak}}(h)\coloneqq\sup_{n=0,1,\ldots,N}\left|\mathbb{E}[\Psi(Y_{n})]-\mathbb{E}[\Psi(y(t_{n}))]\right|,
\end{equation}
where $\Psi:\mathbb{R}^{d}\rightarrow \mathbb{R}$. $\Psi$ is typically taken to be an algebraic polynomial of given degree in view of comparing higher order moments of $y$ (provided they exist, see, e.g., \cite{KloedenPlaten1999}). Two instances of such method taken from \cite{rossler2009} are described by the Butcher's tableaus reported in Tables \ref{RI1table} and \ref{RI6table}. Both schemes are explicit, have $s=3$ stages and are used for the numerical tests in Section \ref{sec:tests} (note that 3 stages is the minimum required to get weak order 2).
\begin{table}[h]
\centering
\caption{Butcher's tableau of the RK method named RI1 in \cite{rossler2009}.}
\label{RI1table}
\centering
\[
\begin{array}{c|rrr|rrr|rrr}
 0 & & & & & & & & & \\[4pt]
 \dfrac{2}{3} & \dfrac{2}{3} & & & 1 & & & & & \\[8pt]
 \dfrac{2}{3} & -\dfrac{1}{3} & 1 & & 0 & 0 & & & & \\[8pt]
\hline
 0 & & & & & & & & & \\[4pt]
 1 & 1 & & & 1 & & & & & \\[4pt]
 1 & 1 & 0 & & -1 & 0 & & & & \\[2pt]
\hline
 0 & & & & & & & & & \\[4pt]
 0 & 0 & & & 1 & & & & & \\[4pt]
 0 & 0 & 0 & & -1 & 0 & & & & \\[2pt]
\hline
\rule{0pt}{15pt} & \dfrac{1}{4} & \dfrac{1}{2} & \dfrac{1}{4} & \dfrac{1}{2} & \dfrac{1}{4} & \dfrac{1}{4} & 0 & \dfrac{1}{2} & -\dfrac{1}{2}\\[8pt]
\hline
\rule{0pt}{15pt} & & & & -\dfrac{1}{2} & \dfrac{1}{4} & \dfrac{1}{4} & 0 & \dfrac{1}{2} & -\dfrac{1}{2}
\end{array}
\]

\end{table}



\begin{table}[H]
\caption{Butcher's tableau of the RK method named RI6 in \cite{rossler2009}.}
\label{RI6table}
\centering
\[
\begin{array}{c|rrr|rrr|rrr}
 0 & & & & & & & & & \\[4pt]
 1 & 1 & & & 1 & & & & & \\[4pt]
 0 & 0 & 0 & & 0 & 0 & & & & \\[2pt]
\hline
 0 & & & & & & & & & \\[4pt]
 1 & 1 & & & 1 & & & & & \\[4pt]
 1 & 1 & 0 & & -1 & 0 & & & & \\[2pt]
\hline
 0 & & & & & & & & & \\[4pt]
 0 & 0 & & & 1 & & & & & \\[4pt]
 0 & 0 & 0 & & -1 & 0 & & & & \\[2pt]
\hline
\rule{0pt}{15pt} & \dfrac{1}{2} & \dfrac{1}{2} & 0 & \dfrac{1}{2} & \dfrac{1}{4} & \dfrac{1}{4} & 0 & \dfrac{1}{2} & -\dfrac{1}{2}\\[8pt]
\hline
\rule{0pt}{15pt} & & & & -\dfrac{1}{2} & \dfrac{1}{4} & \dfrac{1}{4} & 0 & \dfrac{1}{2} & -\dfrac{1}{2}
\end{array}
\]
\end{table}

\section{Extension to SDDEs}
\label{sec:rosslerSDDE}
We now introduce the proposed  extension of the class of methods \cite{rossler2009} summarized in Section \ref{sec:rosslerSODE} to SDDEs defined by discrete delays, i.e., that can be expressed in the form
\begin{equation}\label{SDDE}
\setlength\arraycolsep{0.1em}\begin{array}{rcl}
\dd y(t)&=&f(t,y(t),y(t-\tau_{1}),\ldots,y(t-\tau_{p}))\dd t\\[2mm]
    &&+\displaystyle\sum_{k=1}^{m}g_{k}({\color{black}{t},}y(t),y(t-\tau_{1}),\ldots,y(t-\tau_{p}))\dd W_{k}(t),
    \end{array}
\end{equation}
where $p$ is a positive integer and $0<\tau_{1}<\cdots<\tau_{p}$ are $p$ ordered constant delays. With respect to \eqref{SODE} we have now $f,g_{k}\colon[0,T]\times\mathbb{R}^{d(p+1)}\rightarrow\mathbb{R}^{d}$. We assume throughout that $f$ and $g_{k}$, $k=1,\ldots,m$, are Lipschitz continuous and satisfy a standard linear growth condition with respect to their arguments beyond $t$. This guarantees existence and uniqueness of a strong solution of the IVP for \eqref{SDDE} obtained by assigning a function $[-\tau_{p},0]\ni\theta\mapsto\phi(\theta)$ as initial history, i.e., $y(\theta)=\phi(\theta)$ for $\theta\in[-\tau_{p},0]$, with $\phi(\theta)\in\mathbb{R}^{d}$ a random variable (see, e.g., \cite[Section 5.3]{mao2007stochastic}).

Following the treatment in Section \ref{sec:rosslerSODE} we assume to work with a uniform mesh. To this aim we further assume that the $p$ delays are commensurable with each other and, moreover, that the step size $h$ is such that, for all $q=1,\ldots,p$, the quantity $\kappa_{q}\coloneqq\tau_{q}/h$ is an integer. For simplicity, we also let $N\coloneqq T/h$ be an integer and define again $t_{n}\coloneqq nh$ for $n=0,1,\ldots,N$. Then, by relying on the same random variables introduced in \eqref{Ihat}--\eqref{Ihat2}, the following is a natural extension of R{\"o}{\ss}ler's method \eqref{rosslerSODE:0}--\eqref{rosslerSODE:3} to SDDEs of the form \eqref{SDDE} that reads
\begin{equation}\label{rosslerSDDE:0}
\setlength\arraycolsep{0.1em}\begin{array}{rcl}
Y_{n+1}=Y_{n}&+&\displaystyle h\sum_{i=1}^{s}b_{i}f\left(t_{n}+c_{i}^{(0)}h,H_{n,i}^{(0)},Y_{n-\kappa_{1}+c_{i}^{(0)}},\ldots,Y_{n-\kappa_{p}+c_{i}^{(0)}}\right)\\[4mm]
&+&\displaystyle\sum_{i=1}^{s}\sum_{k=1}^{m}\beta_{i}^{(1)}g_{k}\left(t_{n}+c_{i}^{(1)}h,H_{n,i}^{(k)},Y_{n-\kappa_{1}+c_{i}^{(1)}},\ldots,Y_{n-\kappa_{p}+c_{i}^{(1)}}\right)\hat{I}_{k}\\[4mm]
&+&\displaystyle\sum_{i=1}^{s}\sum_{k=1}^{m}\beta_{i}^{(2)}g_{k}\left(t_{n}+c_{i}^{(1)}h,H_{n,i}^{(k)},Y_{n-\kappa_{1}+c_{i}^{(1)}},\ldots,Y_{n-\kappa_{p}+c_{i}^{(1)}}\right)\frac{\hat{I}_{k,k}}{\sqrt h}\\[4mm]
&+&\displaystyle\sum_{i=1}^{s}\sum_{k=1}^{m}\beta_{i}^{(3)}g_{k}\left(t_{n}+c_{i}^{(2)}h,\hat{H}_{n,i}^{(k)},Y_{n-\kappa_{1}+c_{i}^{(2)}},\ldots,Y_{n-\kappa_{p}+c_{i}^{(2)}}\right)\hat{I}_{k}\\[4mm]
&+&\displaystyle\sum_{i=1}^{s}\sum_{k=1}^{m}\beta_{i}^{(4)}g_{k}\left(t_{n}+c_{i}^{(2)}h,\hat{H}_{n,i}^{(k)},Y_{n-\kappa_{1}+c_{i}^{(2)}},\ldots,Y_{n-\kappa_{p}+c_{i}^{(2)}}\right)\sqrt{h},
\end{array}
\end{equation}
where, for $i=1,\ldots,s$ and $k=1,\ldots,m$,
\begin{equation}\label{rosslerSDDE:1}
\setlength\arraycolsep{0.1em}\begin{array}{rcl}
H_{n,i}^{(0)}=Y_{n}&+&\displaystyle h\sum_{j=1}^{s}A_{i,j}^{(0)}f\left(t_{n}+c_{j}^{(0)}h,H_{n,j}^{(0)},Y_{n-\kappa_{1}+c_{j}^{(0)}},\ldots,Y_{n-\kappa_{p}+c_{j}^{(0)}}\right)\\[4mm]
&+&\displaystyle\sum_{j=1}^{s}\sum_{l=1}^{m}B_{i,j}^{(0)}g_{l}\left(t_{n}+c_{j}^{(1)}h,H_{n,j}^{(l)},Y_{n-\kappa_{1}+c_{j}^{(1)}},\ldots,Y_{n-\kappa_{p}+c_{j}^{(1)}}\right)\hat{I}_{l},
\end{array}
\end{equation}
\begin{equation}\label{rosslerSDDE:2}
\setlength\arraycolsep{0.1em}\begin{array}{rcl}
H_{n,i}^{(k)}=Y_{n}&+&\displaystyle h\sum_{j=1}^{s}A_{i,j}^{(1)}f\left(t_{n}+c_{j}^{(0)}h,H_{n,j}^{(0)},Y_{n-\kappa_{1}+c_{j}^{(0)}},\ldots,Y_{n-\kappa_{p}+c_{j}^{(0)}}\right)\\[4mm]
&+&\displaystyle\sum_{j=1}^{s} B_{i,j}^{(1)}g_{k}\left(t_{n}+c_{j}^{(1)}h,H_{n,j}^{(k)},Y_{n-\kappa_{1}+c_{j}^{(1)}},\ldots,Y_{n-\kappa_{p}+c_{j}^{(1)}}\right)\sqrt{h},
\end{array}
\end{equation}
\begin{equation}\label{rosslerSDDE:3}
\setlength\arraycolsep{0.1em}\begin{array}{rcl}
\hat{H}_{n,i}^{(k)}=Y_{n}&+&\displaystyle h\sum_{j=1}^{s}A_{i,j}^{(2)}f\left(t_{n}+c_{j}^{(0)}h,H_{n,j}^{(0)},Y_{n-\kappa_{1}+c_{j}^{(0)}},\dots,Y_{n-\kappa_{p}+c_{j}^{(0)}}\right)\\[4mm]
&+&\displaystyle\sum_{j=1}^{s}\sum_{\substack{l=1 \\ l\ne k}}^{m}B_{i,j}^{(2)}g_{l}\left(t_{n}+c_{j}^{(1)}h,H_{n,j}^{(l)},Y_{n-\kappa_{1}+c_{j}^{(1)}},\ldots,Y_{n-\kappa_{p}+c_{j}^{(1)}}\right)\frac{\hat{I}_{k,l}}{\sqrt{h}}.
\end{array}
\end{equation}
    Above, for every $q=1,\ldots,p$ and $i=1,\ldots,s$, $Y_{n-\kappa_{q}+c_{i}^{(\ast)}}$ is the approximation of $y(t_{n}-\tau_{q}+c_{i}^{(\ast)}h)$ {\color{black}already computed in a previous step. Indeed, since the delays are commensurable and $\kappa_q$ is an integer for every $q$, the point $t_n -\tau_q+ c_i^{(\ast)}h$ occupies the same relative position within the interval $[t_{n-\kappa_q}, t_{n-\kappa_q+1}]$ as the stage abscissa $t_n + c_i^{(*)}h$ does within the interval $[t_n, t_{n+1}]$ of the current step.  Hence, the quantity $Y_{n-\kappa_{q}+c_{i}^{(0)}}$ denotes the quantity $H_{\ell,i}^{(0)}$ for $\ell$ such that $t_{n}-\tau_q=t_{\ell}$, while the instance of $Y_{n-\kappa_{q}+c_{i}^{(1)}}$ as an argument of $g_k$ denotes the quantity $H_{\ell,i}^{(k)}$ and the instance of $Y_{n-\kappa_{q}+c_{i}^{(2)}}$ as an argument of $g_k$ denotes the quantity $\hat H_{\ell,i}^{(k)}$\footnote{{\color{black}Although $Y_{n-\kappa_{q}+c_{i}^{(\ast)}}$ depends on noise index $k$ for $\ast \in \{1, 2\}$, we omit the superscript $k$ in $Y_{n-\kappa_{q}+c_{i}^{(\ast)}}$ for notational simplicity. }}.} {\color{black}
Regardless of whether
{$Y_{n-\kappa_{q}+c_{i}^{(\ast)}}$ is a mesh value (i.e., $c_{i}^{(\ast)}\in\{0,1\}$) or not, no interpolation of the numerical solution is therefore required.}}

{\color{black}In other words, a step of \eqref{rosslerSDDE:0}-\eqref{rosslerSDDE:3}
has the same form as a step of \eqref{rosslerSODE:0}-\eqref{rosslerSODE:3} that depends on additional but already determined inputs. Moreover, all breaking points belong to the mesh. In the deterministic setting, the absence of interpolation and the inclusion of the breaking points in the mesh are precisely the conditions under which Runge-Kutta methods for DDEs with constant delays retain the order of the underlying ODE method \cite{BellenMasetZennaroGuglielmi2009,BellenZennaro2003}. However, transferring these arguments to the stochastic case is not at all obvious, since the inputs defined by the delayed arguments are not smooth functions but rather Wiener processes, correlated with the present step through shared Brownian motions. Additionally, since individual paths are not smooth, the role of the breaking points in smoothening the expectation of the solution should be analyzed in order to justify the order preservation rigorously. For the purpose of this paper, we therefore regard the preservation of the order as a conjecture, supported by the experiments of Section \ref{sec:tests}. We remark that the schemes reported in Table \ref{RI1table} and Table \ref{RI6table}, used in all experiments therein, fulfill the hypotheses of \cite[Theorem 5.1]{rossler2009}.}
\section{Numerical validation}
\label{sec:tests}
With the aim at providing experimental evidence of the weak order 2 of the method \eqref{rosslerSDDE:0}--\eqref{rosslerSDDE:3}, we illustrate next a number of tests concerning SDDEs of the form \eqref{SDDE}. For each test case we analyze the convergence by evaluating the weak error \eqref{weakerror} with respect to a reference solution $y$ and a uniform mesh for decreasing step size $h=2^{-j}$, $j=0,1,\ldots,5$. In the test case in Subsection~\ref{sec:linear} there exists an analytic expression for the reference solution, and we use that expression for $y$ in \eqref{weakerror}. In the remaining examples, such an explicit expression does not exist, and therefore we approximate $y$ by running $M$ realizations of the numerical method \eqref{rosslerSDDE:0}--\eqref{rosslerSDDE:3} with \eqref{Ihat}--\eqref{Ihat2} using a much smaller step size $h_{\text{ex}}$, ranging from $2^{-10}$ to $2^{-16}$ depending on the test case. For each test case and value of $h$ considered, $M$ also indicates the number of realizations of the corresponding numerical solution. Moreover, in every test and without loss of generality we restrict to measure the weak error at the final time $T$, hence
\begin{equation}\label{weakerrorend}
E_{\text{weak}}(h)\coloneqq\left|\mathbb{E}[\Psi(Y_{N})]-\mathbb{E}[\Psi(y(t_{N}))]\right|.
\end{equation}

In the following tests we consider a rather large number $M$ of realizations (ranging from $10^{7}$ to $4.5\times 10^{7}$) in order to obtain a clear confirmation of the weak order 2. To this aim, MATLAB codes are implemented separately for each specific SDDE to improv{\color{black}e} computational efficiency and thus enabling us to obtain such large numbers of realizations. All relevant scripts are available at \url{https://cdlab.uniud.it/software#w_ord_2_RK_SDDEs}, implemented in MATLAB R2025b. All relevant tests are r{\color{black}u}n on a PC with the following specifications: Intel(R) Core(TM) Ultra 7 256V, RAM 16GB and Graphics card 128 MB. Finally, in Section \ref{sec:linear} we use both the methods described in Table \ref{RI1table} and Table \ref{RI6table}. As the results are the same (for this and other tests that we omit to report here), in the remaining experiments in Section \ref{sec:nonlinear}, Section \ref{sec:system} and Section \ref{sec:supply} we restrict to the method in Table \ref{RI6table} without loss of generality.

\subsection{A linear SDDE}
\label{sec:linear}
The first example we consider is taken from \cite{kuchlerplaten} and is the scalar linear SDDE with a single multiplicative noise and a single constant delay
\begin{equation}\label{eq:linear}
\dd y(t)=(\mu_{1}y(t)+\mu_{2}y(t-1))\dd t+\sigma y(t)\dd W(t)
\end{equation}
with $\mu_{1},\mu_{2},\sigma\in\mathbb{R}$. With respect to \eqref{SDDE}
we have $d=m=p=1$ and, in particular, $\tau_{1}=1$. The  IVP of interest is obtained by setting $y(t)=y(0)$ for $t\in[-1,0]$. Correspondingly, an exact\footnote{Note that in order to obtain a reference solution for the weak error \eqref{weakerrorend} we use a uniform left-hand rectangles quadrature of the integral in \eqref{kuchlersol1} with step size $h_{\text{ex}}$.\label{foonote1}} expression for the solution can be described recursively for $k=1,2,\ldots$ as
\begin{equation}\label{kuchlersol1}
y(t)=\Phi_{t,k-1}\left(y(k-1)+\int_{k-1}^{t}\mu_{2}y(s-1)\Phi^{-1}_{s,k-1}\dd s\right),\quad t\in[k-1,k],
\end{equation}
where
\begin{equation}\label{kuchlersol2}
\Phi_{t,t_{0}}\coloneqq\text{exp}\left(\left(\mu_{1}-\frac{1}{2}\sigma^2 \right)(t-t_{0})+\sigma(W(t)-W(t_{0}))\right),\quad t\in[t_{0},t_{0}+1].
\end{equation}
     
\begin{figure}[H]
    \centering
    \begin{minipage}{0.49\textwidth}
        \centering
        \includegraphics[width=\linewidth, height=0.25\textheight]{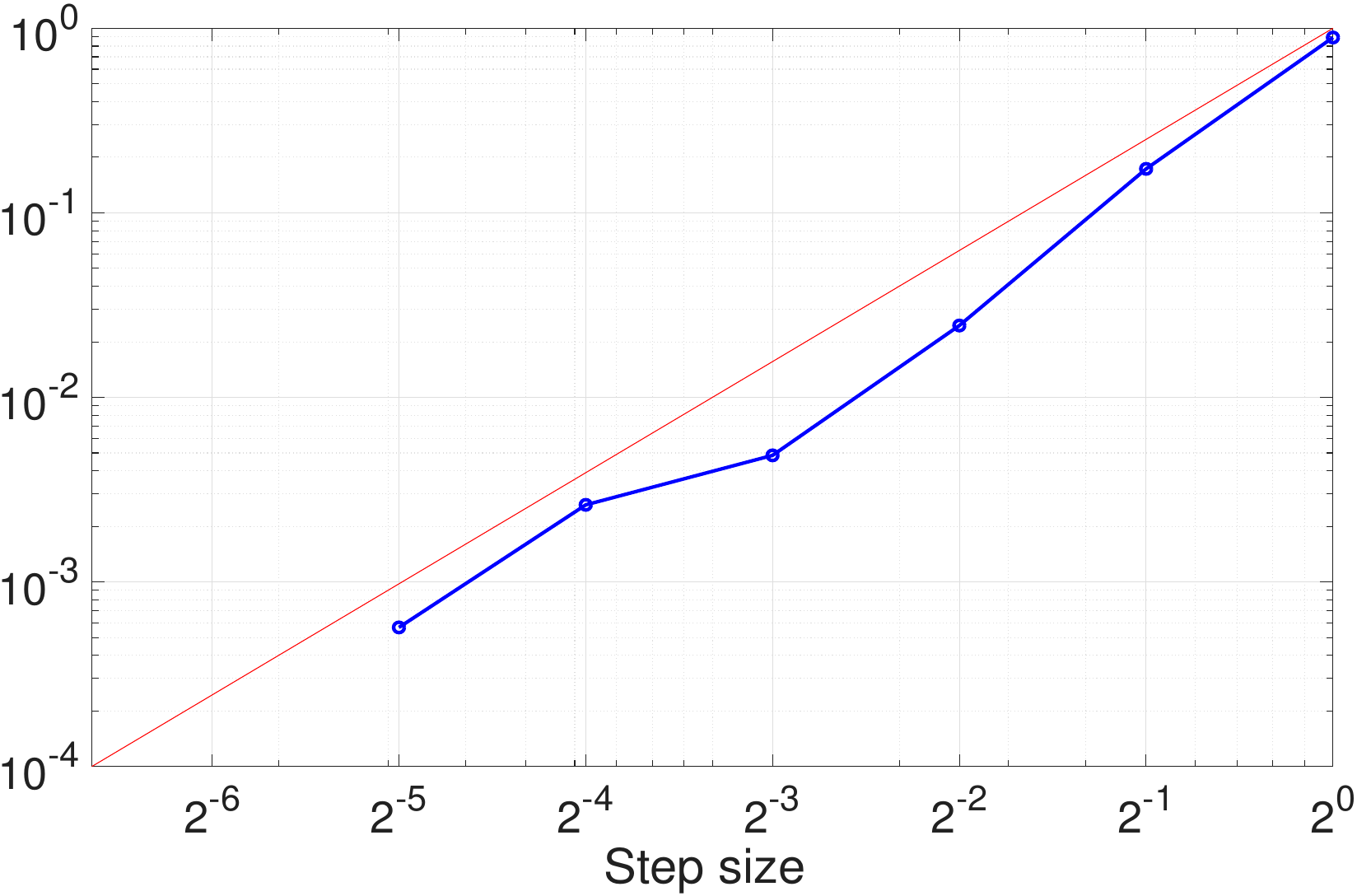}
    \end{minipage}
    \hfill
    \begin{minipage}{0.49\textwidth}
        \centering
        \includegraphics[width=\linewidth, height=0.25\textheight]{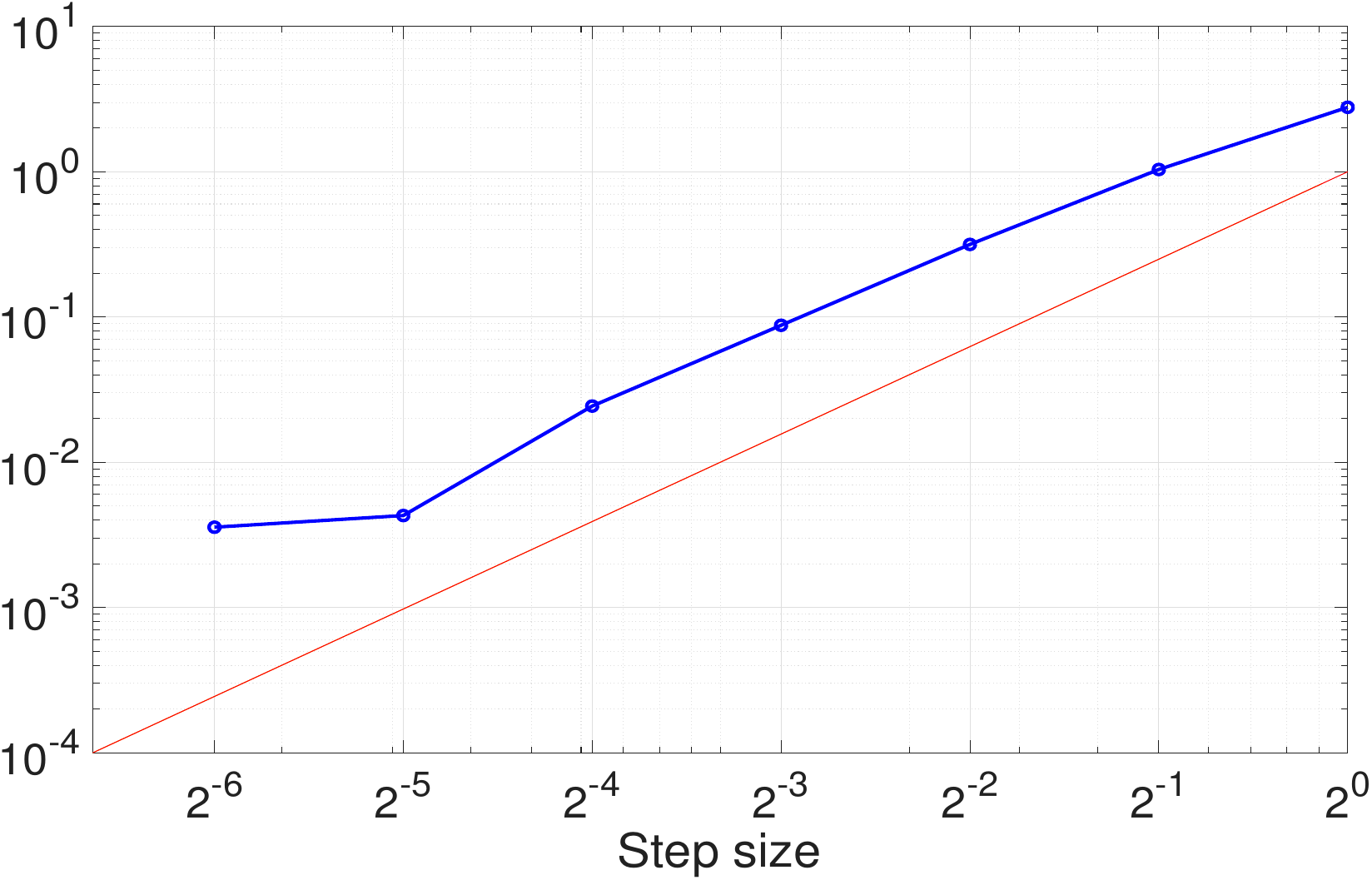}
    \end{minipage}
    
    \caption{weak order analysis of the extended R{\"o}{\ss}ler methods of Table \ref{RI1table} (left) and Table \ref{RI6table} (right) for the linear SDDE \eqref{eq:linear} with $\mu_{1}=\mu_{2}=1$, $\sigma=0.5$, $y(0)=1$, $T=2$, $M=2\times 10^{7}$ and $h_{\text{ex}}=2^{-14}$ (thin red: reference line of slope 2).}   
    \label{fig:linear1}
\end{figure}

The results of a first test are reported in Figure \ref{fig:linear1} about comparing the convergence of the weak error of the extended R{\"o}{\ss}ler methods of Table \ref{RI1table} (left panel) and Table \ref{RI6table} (right panel): weak order 2 is confirmed for both schemes. The experiment refers to the the model parameter values $\mu_{1}=\mu_{2}=1$ and $\sigma=0.5$, initial value $y(0)=1$ and $T=2$. The weak endpoint error \eqref{weakerrorend} is calculated with $M=2\times 10^{7}$ realizations and the reference solution is obtained from \eqref{kuchlersol1} with \eqref{kuchlersol2}  by using a left-hand rectangles rule for the integral in \eqref{kuchlersol1} with a step size $h_{\text{ex}}=2^{-14}$ (recall footnote \ref{foonote1}).
\vspace{0.4 cm}

\begin{figure}[H]
    \centering
    \begin{minipage}{0.49\textwidth}
        \centering
        \includegraphics[width=\linewidth, height=0.25\textheight]{LatestNEWRI1-Kuchler-T-2-sigma-0point5-15000000.pdf}
    \end{minipage}
    \hfill
    \begin{minipage}{0.49\textwidth}
        \centering
        \includegraphics[width=\linewidth, height=0.25\textheight]{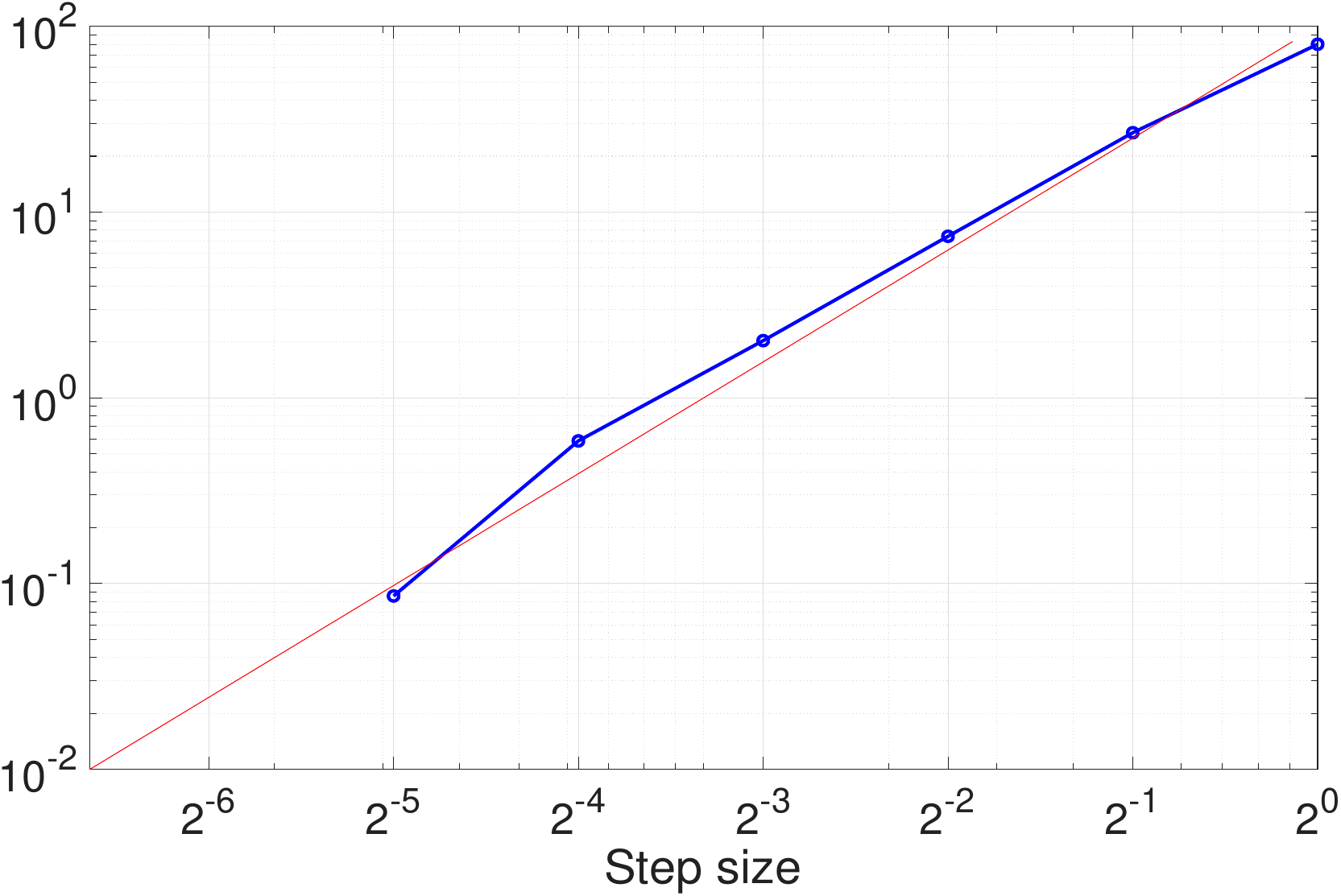}
    \end{minipage}
    \caption{weak order analysis of the extended R{\"o}{\ss}ler method of Table \ref{RI6table} with respect to \eqref{weakerrorend} with $\Psi(y)=\Psi_{1}(y)\coloneqq y$ (left) and $\Psi(y)=\Psi_{2}(y)\coloneqq y^{2}$ (right) for the li{\color{black}ne}ar SDDE \eqref{eq:linear} with $\mu_{1}=\mu_{2}=1$, $\sigma=0.5$, $y(0)=1$, $T=2$, $M=2\times 10^{7}$ and $h_{\text{ex}}=2^{-14}$ (thin red: reference line of slope 2).}   
    \label{fig:linear2}
\end{figure}

The results of a second similar test are reported in Figure \ref{fig:linear2} for the same model and numerical parameters. Here only the extended R{\"o}{\ss}ler method in Table \ref{RI6table} is used but \eqref{weakerrorend} refers to either $\Psi(y)=\Psi_{1}(y)\coloneqq y$ (left panel) and $\Psi(y)=\Psi_{2}(y)\coloneqq y^{2}$ (right panel). It is clear from the results that the weak order of convergence remains the same while the error constant changes with $\Psi$ \cite{KloedenPlaten1999}.
{\color{black} We conclude this section with a brief discussion on statistical errors and the computational effort with respect to increasing number $m$ of noise terms.
We conducted an experiment with the linear SDDE \eqref{eq:linear} using $M=8 \times 10^{7}$ realizations and $h_{\text{ex}}=2^{-14}$ to obtain the confidence intervals for the weak errors for each step size, reported in Table~\ref{tab:placeholder}. We observe that, for all but the smallest step size considered, the confidence intervals are significantly smaller than the weak errors themselves. For $h=2^{-5}$, however, the confidence interval becomes comparable to the estimated weak error, implying that the error estimate at this step size should be interpreted with some caution. Nevertheless, the convergence behavior observed for the remaining step sizes already provides a reliable indication of the convergence order, and the result with the smallest one is consistent with this behavior, although its precision is more limited. Similar analyses can be conducted for all other upcoming models.
\begin{table}[htbp]
\end{table}
\begin{table}[H]
    \centering
    \renewcommand{\arraystretch}{1.3}
     \caption{weak errors and corresponding confidence intervals.}
    \begin{tabular}{l|c|r}
         \textbf{Step Size}&\textbf{Weak Error} & \textbf{Confidence Interval}  \\ \hline
         $2^{-5}$ &  0.0056&  [0.0025,0.0087]\\ \hline
         $2^{-4}$ &0.0254 &[0.0223,0.0285]\\ \hline
         $2^{-3}$ &  0.0902& [0.0871,0.0933] \\ \hline
         $2^{-2}$ &0.3213 & [0.3183,0.3243]\\ \hline
         $2^{-1}$ &1.0357 & [1.0328,1.0385]\\ \hline
         $2^{0}$ & 2.7780 & [2.7754,2.7806]\\ \hline
         
    \end{tabular}
   
    \label{tab:placeholder}
\end{table}

In addition to the statistical errors we measure the effect of increasing noise terms on the computational cost of the method. We replaced $\sigma y(t)\dd W(t)$ in \eqref{eq:linear} with $\Sigma_{i=1}^{m}\sigma y(t) \dd W_{i}$ for $m=1,2, \dots, 5$ and measure the computational cost for $h_{\text{ex}}=2^{-14}$ and $M=2 \times 10^{6}$, Figure \ref{f3m}: the cost is roughly proportional to $m$, consistently with the theory in \cite{rossler2009} (recall footnote \ref{noisefootnote}).
\begin{figure}[H]
    \centering
    \includegraphics[width=0.5\linewidth, height=0.25\textheight]{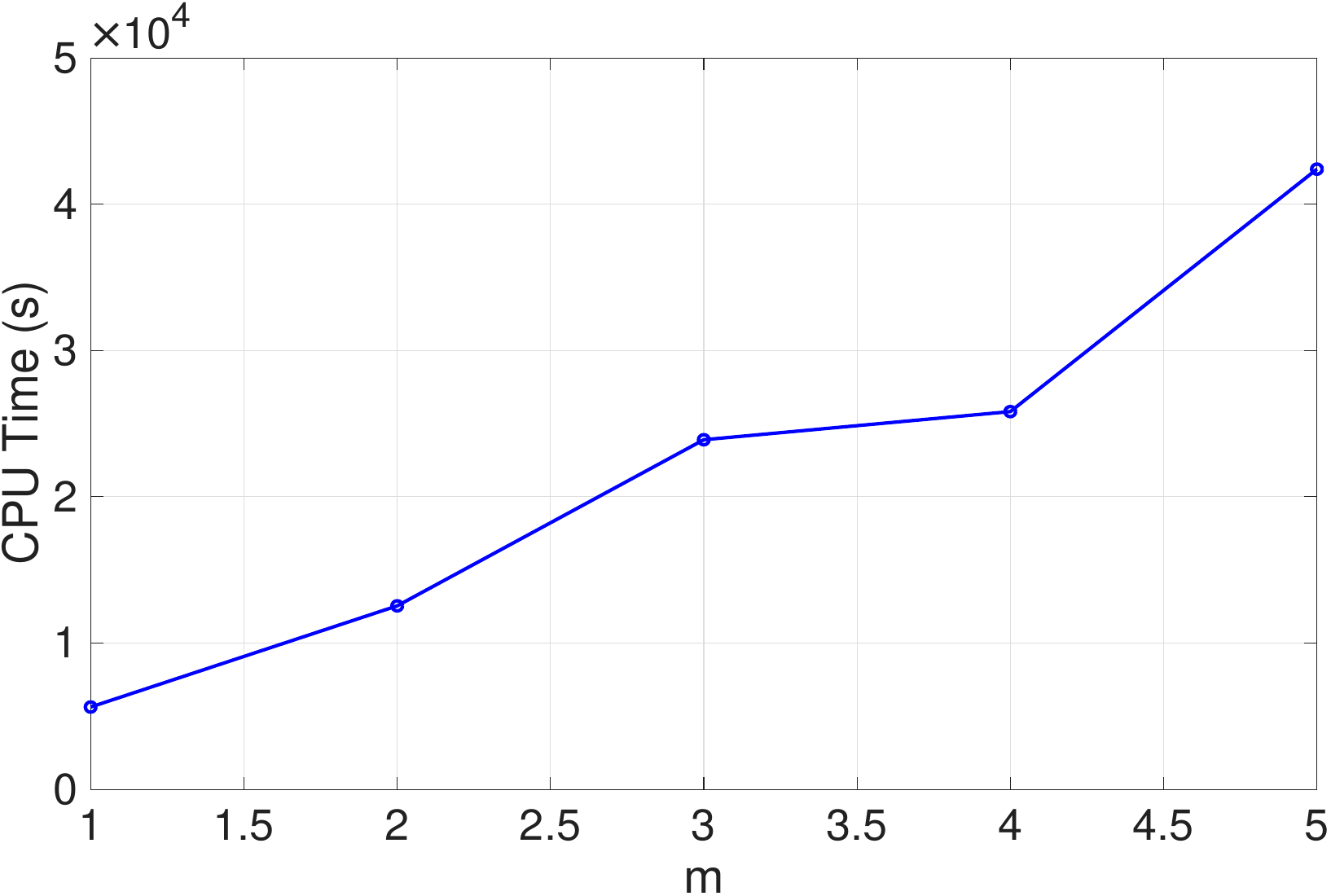}
\caption{CPU time in seconds plotted against the number of noise terms.}\label{f3m}  
\end{figure}
} 
\subsection{A nonlinear SDDE}
\label{sec:nonlinear}
The second example we consider is a logistic model inspired by \cite{Wu2012StochasticDelayLogistic} and is the scalar SDDE with nonlinear drift, a single multiplicative noise and a single constant delay
\begin{equation}\label{eq:nonlinear}
\dd y(t)=\mu y(t)(1-y(t-1))\dd t+\sigma y(t)\dd W(t)
\end{equation}
with {\color {red}$\mu,\sigma\in\mathbb{R}$.} With respect to \eqref{SDDE}
we have again $d=m=p=1$ and $\tau_{1}=1$. A corresponding IVP is obtained by assigning an initial history $\phi(t)\in\mathbb{R}$ for $t\in[-1,0]$.

\begin{figure}[H]
    \centering
    \includegraphics[width=0.5\linewidth, height=0.25\textheight]{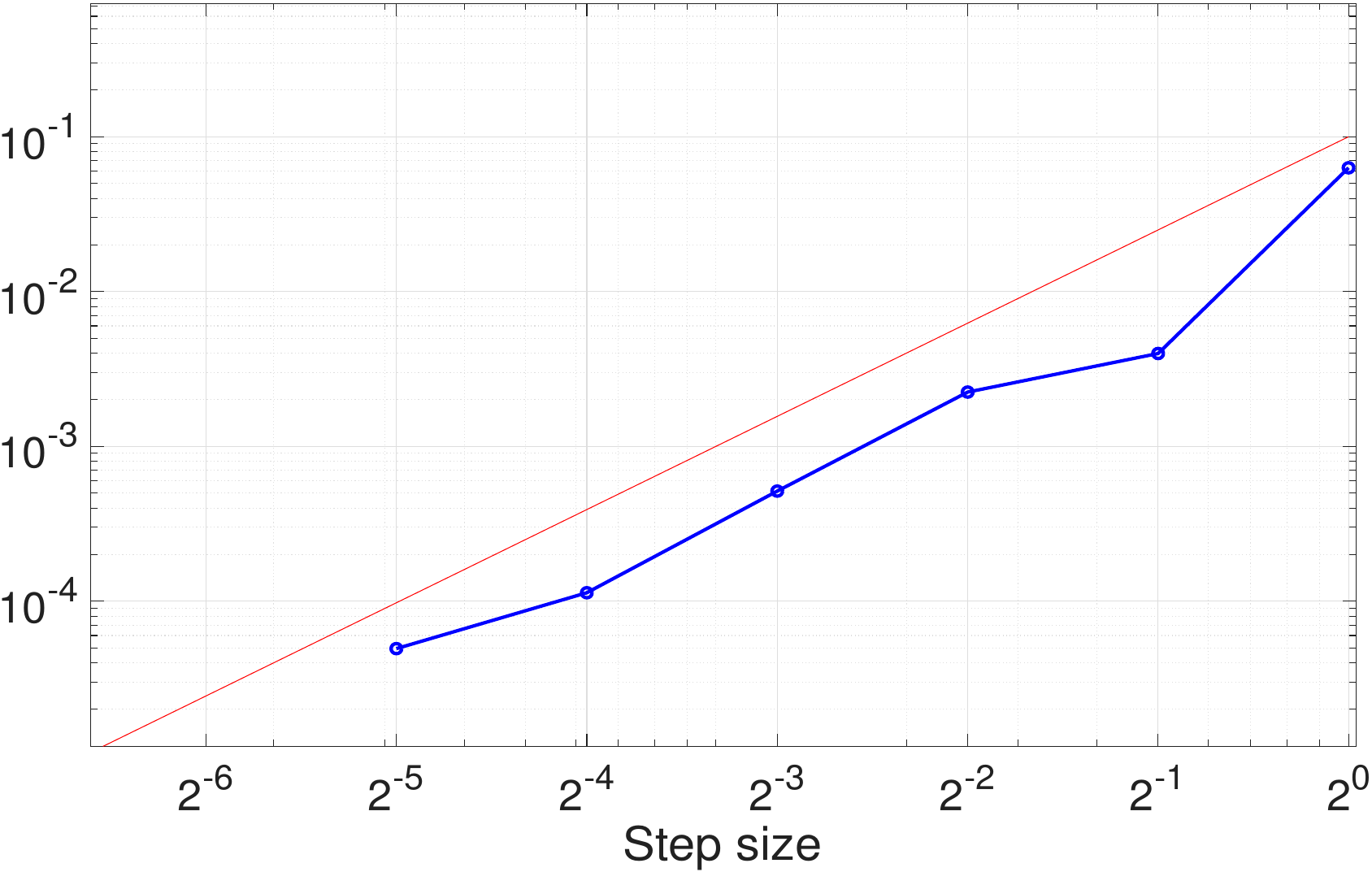}
    \caption{weak order analysis of the extended R{\"o}{\ss}ler method of Table \ref{RI6table} for the logistic SDDE \eqref{eq:nonlinear} with $\mu=1$, $\sigma=0.5$, $\phi(t)=\cos(t)$ for $t \in [-1,0]$, $T=2$, $M=4.5\times 10^{7}$ and $h_{\text{ex}}=2^{-16}$ (thin red: reference line of slope 2).}\label{fig:nonlinear}
\end{figure}

The results in Figure \ref{fig:nonlinear} concern the extended R{\"o}{\ss}ler method of Table \ref{RI6table} and confirm weak order 2 for the weak endpoint error \eqref{weakerrorend}. The experiment refers to the model parameter values $\mu=1$ and $\sigma=0.5$, initial history $\phi(t)=\cos(t)$ and $T=2$. The error is calculated with $M=4.5\times 10^{7}$ realizations and the reference solution is obtained with the same scheme with step size $h_{\text{ex}}=2^{-16}$.

\subsection{A system of SDDEs with multiple noises}
\label{sec:system}
The third example we consider is the predator-prey model from  \cite{predatorprey2015} and is the system of SDDEs with nonlinear drift, two multiplicative noises and a single constant delay
\begin{equation}\label{eq:system}
\left\{\setlength\arraycolsep{0.1em}\begin{array}{rcl}
\dd y_{1}(t)&=&y_{1}(t)[\alpha-\beta y_{1}(t)-\gamma y_{2}(t-\tau)]\dd t+\sigma_{1} y_{1}(t)\dd W_{1}(t)\\[2mm]
\dd y_{2}(t)&=&y_{2}(t)[-\delta+\kappa y_{1}(t-\tau)]\dd t+\sigma_{2}y_{2}(t)\dd W_{2}(t),
\end{array}
\right.
\end{equation}  
where $y_{1}$ and $y_{2}$ represent the densities of prey and predator respectively, $W_{1}$ and $W_{2}$ are independent scalar Brownian motions and $\alpha,\beta,\gamma,\delta,\sigma_{1},\sigma_{2}\in\mathbb{R}$. With respect to \eqref{SDDE} we have now $d=m=2$, $p=1$ and $\tau_{1}=\tau$. A corresponding IVP is obtained by assigning an initial history $\phi(t)\in\mathbb{R}^{2}$ for $t\in[-\tau,0]$.

\begin{figure}[H]
    \centering
    \includegraphics[width=0.5\linewidth, height=0.25\textheight]{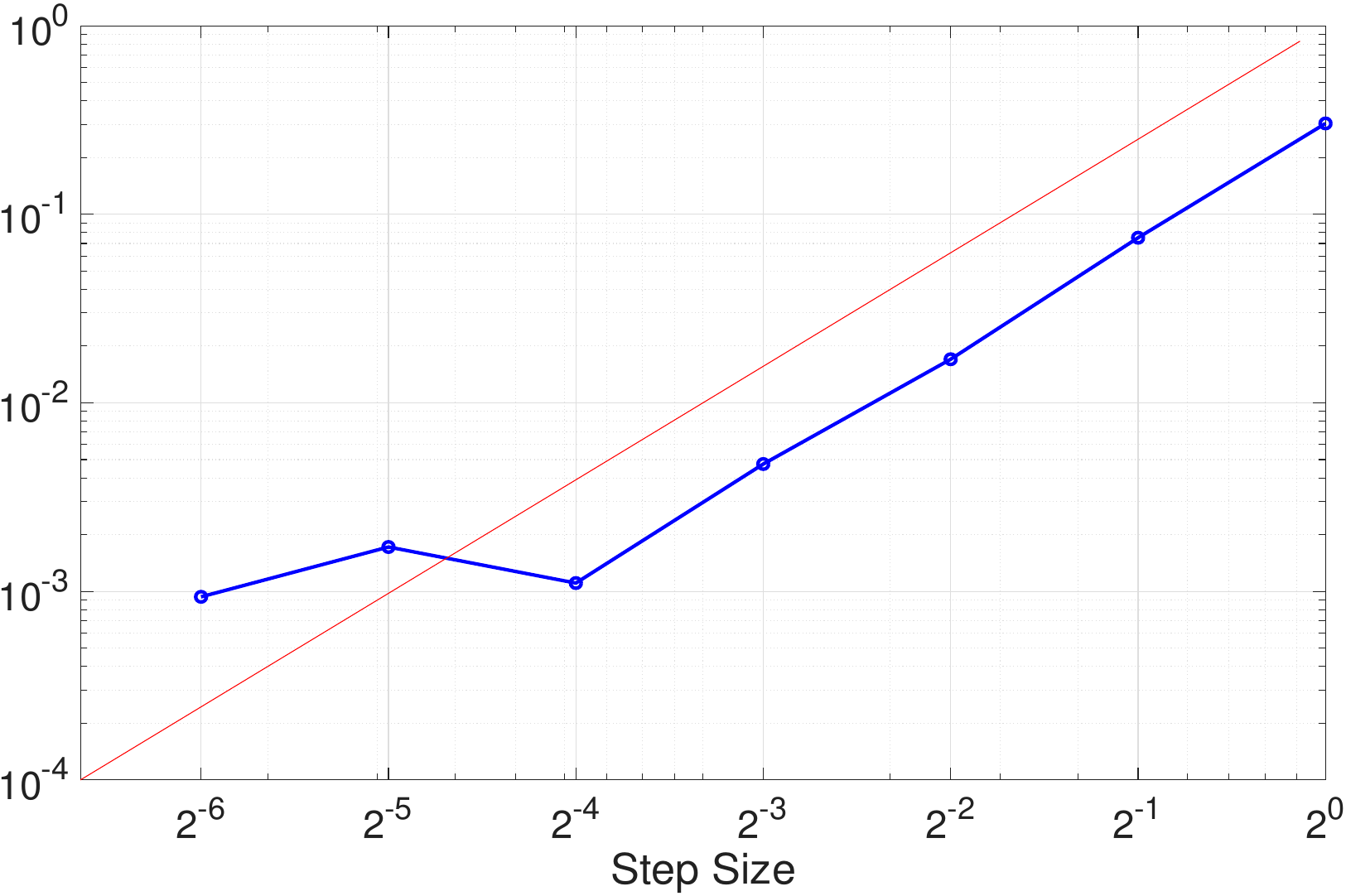}
    \caption{weak order analysis of the extended R{\"o}{\ss}ler method of Table \ref{RI6table} for the predator-prey SDDE \eqref{eq:system} $\alpha=1$, $\delta=0.5$, $\beta=\gamma=\kappa=0.1$, $\sigma_{1}=\sigma_{2}=0.5$, $\phi(t)=(5,2)^{T}$ for $t \in [-1,0]$, $T=2$, $M=5\times 10^{6}$ and $h_{\text{ex}}=2^{-10}$ (thin red: reference line of slope 2).}\label{fig:system}
\end{figure}

The results in Figure \ref{fig:system} concern the extended R{\"o}{\ss}ler method of Table \ref{RI6table} and confirm weak order 2 for the weak endpoint error \eqref{weakerrorend}. The experiment refers to the model parameter values $\alpha=1$, $\delta=0.5$, $\beta=\gamma=\kappa=0.1$, $\sigma_{1}=\sigma_{2}=0.5$, initial history $\phi(t)=(5,2)^{T}$ for $t \in [-1,0]$ and $T=2$. The error is calculated with {\color{black}$M=5 \times  10^{6}$} realizations and the reference solution is obtained with the same scheme with step size {\color{black}$h_{\text{ex}}=2^{-10}$.}

\subsection{A supply-chain model with multiple delays}
\label{sec:supply}
The last example we consider is a supply-chain model adapted from \cite{supplychain2020} describing the evolution of
a brand goodwill, and is the scalar linear SDDE with multiple delays
\begin{equation}\label{eq:supply}
\dd y(t)=(\alpha y(t-\tau_{M})+\beta y(t-\tau_{R})-\gamma y(t-\tau))\dd t+\delta y(t-\tau)\dd W(t).
\end{equation}
Above $\alpha,\beta,\gamma,\delta\in\mathbb{R}$, while $\tau_{M}$ and $\tau_{R}$ represent, respectively, the manufacturer’s and retailer’s additional advertising delays beyond an unavoidable delay $\tau$. With respect to \eqref{SDDE} we have $d=m=1$ and $p=3$. A corresponding IVP is obtained by assigning an initial history $\phi(t)\in\mathbb{R}$ for $t\in[-\tau_{3},0]$.

\begin{figure}[H]
    \centering
    \includegraphics[width=0.5\linewidth, height=0.25\textheight]{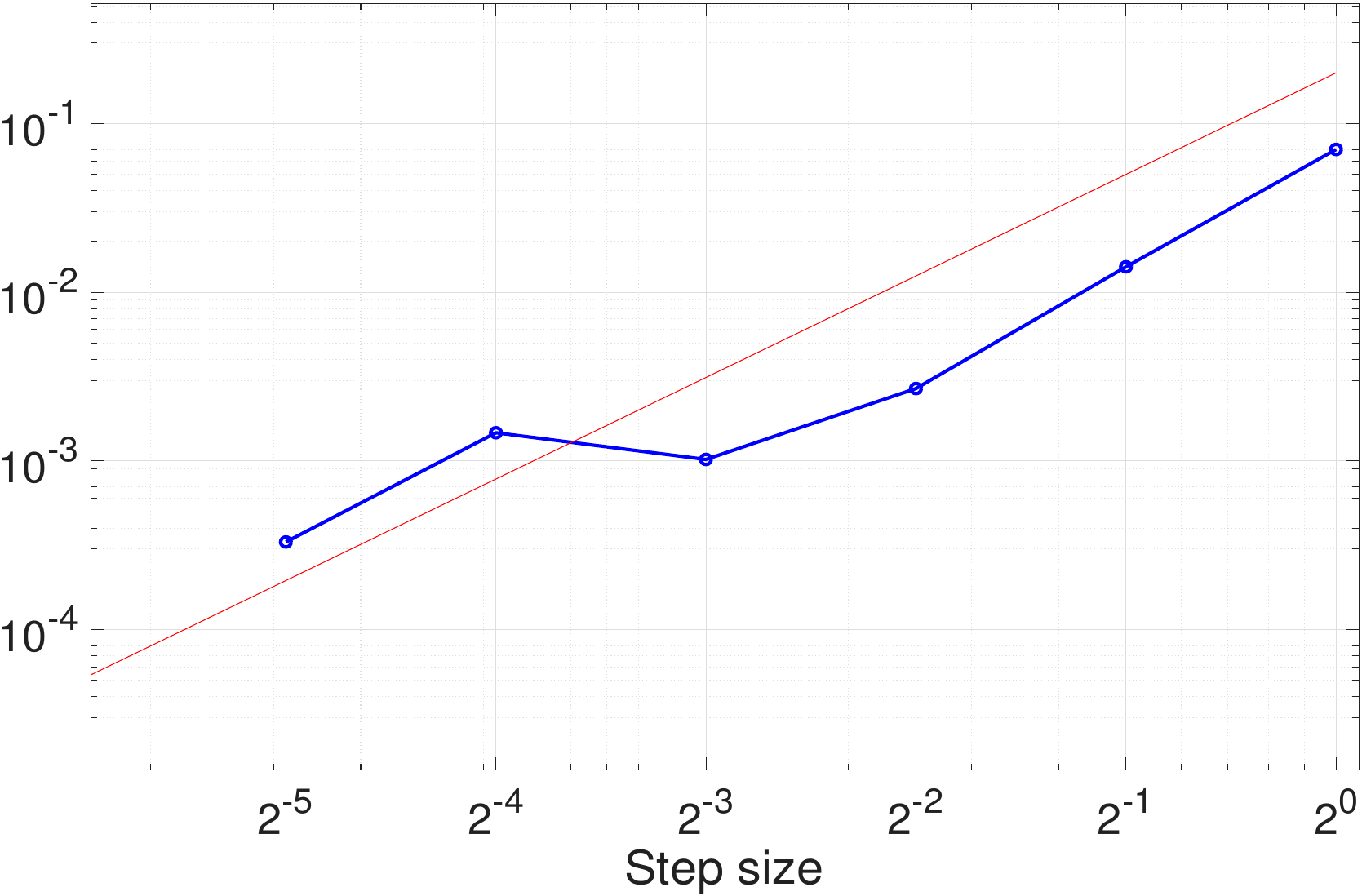}
    \caption{Weak error analysis of \eqref{eq:supply} with $\alpha=\beta=1$, $\gamma=0.8$, $\delta=0.5$, delays $\tau_{M}=5$, $\tau_{R}=4$ and $\tau=1$, initial history $\phi(t)=5$ for $t \in [-5,0]$, $T=5$, $M=1.5\times 10^{7}$ and $h_{\text{ex}}=2^{-10}$ (thin red: reference line of slope 2).}\label{fig:supply}
\end{figure}

The results in Figure \ref{fig:supply} concern the extended R{\"o}{\ss}ler method of Table \ref{RI6table} and confirm weak order 2 for the weak endpoint error \eqref{weakerrorend}. The experiment refers to the model parameter values $\alpha=\beta=1$, $\gamma=0.8$, $\delta=0.5$, delays $\tau_{M}=5$, $\tau_{R}=4$ and $\tau=1$, initial history $\phi_(t)=5$ for $t \in [-5,0]$ and $T=5$. The error is calculated with $M=1.5\times 10^{7}$ realizations and the reference solution is obtained with the same scheme with step size $h_{\text{ex}}=2^{-10}$.
In addition to the above experiment, several other tests were conducted on the same model varying step sizes and number of realizations. A key emerging observation is that the number of realizations must be increased as the step size becomes smaller, confirming what indicated in \cite{higham2021introduction}.

\section{Conclusions and future directions}
\label{sec:conclusions}
In this work we proposed an extension from SODEs to SDDEs of the class of RK methods introduced by R{\"o}{\ss}ler in \cite{rossler2009}, with the aim at constructing derivative-free schemes of weak order 2. This first extension holds for SDDEs with a finite number of constant, discrete and commensurate delays, and works on a uniform mesh with step size a submultiple of each delay. A thorough numerical investigation confirmed the preservation of the weak order 2 for several classe{\color{black}s} of SDDEs, including instances with multiple delays and multiple noises. Relevant MATLAB codes are freely available at \url{https://cdlab.uniud.it/software#w_ord_2_RK_SDDEs}.

As anticipated in the Introduction, this represents our first step towards extending to SDDEs the class of FRK methods \cite{MasetDDE,MasetTorelliVermiglio} developed for deterministic DDEs. Two developments are therefore envisaged in this direction. First, constructing a FRK scheme that includes a treatment of the stochastic part in the way of R{\"o}{\ss}ler and experimentally verifying the preservation of the weak order 2. In these regards the so-called {\it improved Euler method} or {\it explicit trapezoidal rule} from the family of Heun methods described in \cite{MasetDDE} is a natural FRK candidate. Second, a rigorous proof of convergence in the weak sense is necessary to fully grasp the interplay between delay(s) and stochasticity{\color{black}, as discussed in Section \ref{sec:rosslerSDDE}. In particular, it should clarify whether the order conditions of \cite{rossler2009} remain sufficient when the expansions involve delayed arguments, and how the regularity of the expectation of the solution behaves at breaking points}. To note that the theoretical approaches to the convergence analysis adopted in \cite{rossler2009} and in \cite{MasetDDE} are different, the former being based on the development of so-called colored rooted trees \cite{Rossler2004,Rossler2006RootedTreeSRK}, the latter stemming from the so-called Albrecht approach \cite{Albrecht1987}. It is therefore reasonable to expect that a trivial combination of the two strategies is either impossible or not quite so trivial.

In addition to moving forward towards our quest to extend the FRK methods to SDDEs, we also intend to investigate a similar extension of R{\"o}{\ss}ler's work in terms of the strong error, see \cite{Rossler2010}.

\section*{Acknowledgments}
The authors are members of INdAM research group GNCS and UMI; AA and DB are members of UMI research group “Modellistica socio-epidemiologica”. The work of AA and DB was partially supported by the Italian Ministry of University and Research (MUR) through the PRIN 2022 project (No. 20229P2HEA) “Stochastic numerical modelling for sustainable innovation”, Unit of Udine (CUP G53C24000710006). The work of FW was supported by the Italian Ministry of University and Research (MUR) through a PhD grant PNRR DM629/24 (CUP: G23C24001350003).

\bibliographystyle{abbrv}
\bibliography{bibliografia}

\end{document}